\documentclass[11pt]{article}
\usepackage{pgfplots}
\pgfplotsset{compat=1.8}
\usepackage{wrapfig}
\usepackage[utf8]{inputenc}
\usepackage[T1]{fontenc}
\usepackage[english]{babel}
\usepackage{amsmath,amssymb,amsfonts}
\usepackage{amsthm}
\usepackage{geometry}
\usepackage{mathtools}
\usepackage{graphicx}
\mathtoolsset{showonlyrefs}
\usepackage{bm}
\usepackage{subcaption}
\usepackage{algorithm}
\usepackage{algpseudocode}
\usepackage[hidelinks]{hyperref}
\usepackage[shortlabels]{enumitem}
\usepackage{listings}
\usepackage{cite}
\usepackage{titlesec}
\usepackage[draft, todonotes={textsize=scriptsize}]{changes}
\setuptodonotes{color=red, backgroundcolor=white,
  bordercolor=red, tickmarkheight=10pt}
\titleformat{\subsubsection}[runin]{\normalfont\normalsize\itshape}{\thesubsubsection}{5pt}{}

\newtheorem{theorem}{Theorem}

\newtheorem{remark}[theorem]{Remark}

\theoremstyle{definition}

\newcommand{\R}{\mathbb{R}}

\newcommand{\NO}{\mathcal{N}}
\newcommand{\dd}{\mathrm{d}}
\renewcommand{\P}{\mathcal{P}}

\mathtoolsset{showonlyrefs}

\DeclareMathOperator*{\argmin}{arg\,min}

\begin{document}

\title {A Study of Particle Motion in the Presence of Clusters}
\author{Christian Wald\footnotemark[1] 
		\and Gabriele Steidl\footnote{Institute of Mathematics,
	   Technische Universität Berlin,
	   Stra{\ss}e des 17.\ Juni 136, 
	   10623 Berlin, Germany,
	   {\ttfamily\{wald, steidl\}@math.tu-berlin.de},
      \url{http://tu.berlin/imageanalysis}.
	} 
 }
\date{\today}

\maketitle

\begin{abstract}
The motivation for this study came from the task of analysing
the kinetic behavior of single molecules in a living cell
based on Single Molecule Localization Microscopy.
Given measurements of both the motion of clusters and molecules,
the main task consists in detecting if a molecule belongs to a cluster. While the exact size of the clusters is usually unknown, 
upper bounds are available. In this study, we simulate the cluster movement by a Brownian motion 
and those of the particles by a Gaussian mixture model with two modes
depending on the position of the particle within or outside a cluster. We propose various variational models 
to detect if a particle lies within a cluster
based on the Wasserstein and maximum mean discrepancy distances between measures.
We compare the performance of the proposed models
for simulated data.
\end{abstract}

\section{Introduction}
In this paper, we are interested in a setting where particles 
and particle clusters are moving with different ,,speed'' in space, 
which we assume to be $\R^d$. In practice, $d=2$ or $d=3$ will be of interest.
A single particle might enter a cluster and will then be carried by the cluster, 
but is also allowed to leave it at some point. 
A typical real-world setting is 
Single Molecule Localization Microscopy (SMLM) \cite{lelek2021single} 
which includes high precision in localization of both single molecules and clusters.
However, the size of the clusters is usually unknown and only upper bounds exist.
The original motivation to deal with this topic is related to  
the study the consistency and morphology of protein microclusters 
formed in the context of so-called $T$ cell activation, 
which have been a topic of interest in biology over the last 15 years \cite{sherman2013super,rossboth2018tcrs}. 
In particular, it is  not possible to see precisely whether a single molecule is part of a cluster. 

Assuming that the movement patterns of single particles and cluster follow a Brownian motion with quite different parameters, 
we propose a simple simulation of the particle-cluster behavior
using a Gaussian mixture model (GMM) with two modes which address the cluster membership of a particle.
Knowing particle trajectories the parameters of the GMM can be estimated 
by the expectation-minimization (EM) algorithm \cite{dempster1977maximum}.

Based on this model, our main task is to determine whether a particle belongs to some cluster.
We will use a variational approach to determine the probability 
that a particle belongs to a cluster at a certain time.
More precisely, we will define different functionals which minimizers
serve as estimators of the above probability.
These functionals consists of two summands, 
which take i) the distance of the single particle to the closest cluster 
and 
ii) the movement patterns from the GMM into account. 
For the second part, we present several choices:
\begin{itemize}
\item[-] we relate path trajectories of particles
to empirical probability measures
and use distances between probability measures such as the Wasserstein distance and the maximum mean discrepancy (MMD)
to consider their distance to the cluster component of the GMM.
\item[-] we relate the associated empirical probability measures of particle trajectories and cluster trajectories. 
\end{itemize}
We provide numerical examples on simulated data for both approaches and different distances between measures.
Here we recognized that the latter approach obtains a better accuracy when comparing with the ground truth of the simulated data.
\\[1ex]
\textbf{Outline of the paper.}	In Section \ref{sec:model}, we explain our model for simulating particle-cluster motion
and recall the EM algorithm for determining the parameters of the associated two mode GMM if samples of particle trajectories are available.
Then, in Section \ref{sec:clusters}, we introduce the functionals for determining the probability that a particle belongs to a cluster at a given time.
In Section \ref{sec:num} we compare our different approaches numerically for simulated data.
Conclusions are drawn in Section \ref{sec:conc}.

\section{Simulation of single particle versus clusters motion} \label{sec:model}
%
We start by describing our particle-cluster movement simulation in $\mathbb R^d$, 
where we will rely on $d=2$ in our numerical experiments.
Then we will see that the motion of the particle follows a GMM with two modes 
depending on its position inside or outside some cluster.
Given (simulated) measurements of the particle positions at several times, the parameters of this GMM
can be estimated by the EM algorithm.

In the following, we use $\| \cdot\|$ for the Euclidean norm on $\mathbb R^d$
and the abbreviation $[N]:= \{1,\ldots,N\}$.

\subsection{Simulation}
\begin{figure}
\begin{minipage}{0.6\linewidth}

\includegraphics[width=\linewidth]{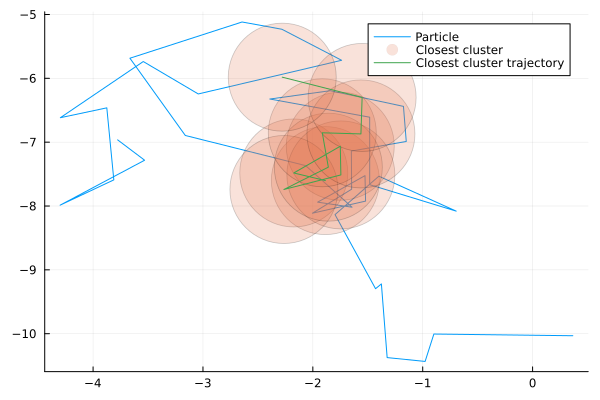}
\end{minipage}
\begin{minipage}[b]{0.4\linewidth}
\caption{Trajectory of a particle together with its closest cluster. The cluster is only displayed if the particle is inside of it. The movement outside and inside the cluster clearly differs.}\label{sim:movment}
\end{minipage}

\end{figure}

We consider $N$ clusters at time points $t \in [T]$
modeled as balls with the same radius $r>0$ and centers $c_{n,t}$,
$n \in [N]$, $t \in [T]$.
The clusters move according to a Brownian motion, i.e., 
\[
\Delta c_n(t) \coloneqq c_{n,t+1}-c_{n,t} \sim \NO(m_c,\sigma^2_c).
\] 
where we write shorthand $\NO(m_c,\sigma^2_c)$ for $\NO(m_c,\sigma^2_c \, I_d)$.
The movement of a single particle is described by its position $p_t\in \R^d$, $t \in [T]$. 
Outside clusters we assume a Brownian motion 
$$\Delta p(t) \coloneqq p_{t+1}-p_{t} \sim \NO(m_{out},\sigma^2_{out} ),$$
where $\sigma_{out} > \sigma_{c}$.
Inside a cluster, the particle is carried by the cluster with an additional small Brownian motion $\NO(0,\sigma^2_{pc} )$. 
This results in an overall Brownian motion within a cluster given by
$$
\Delta p(t)  \sim \NO(m_{in},\sigma^2_{in} ), 
\qquad m_{in} \coloneqq m_c, \, \sigma^2_{in} \coloneqq \sigma^2_c + \sigma^2_{pc}.
$$ 
Note that due to the small Brownian motion inside the cluster, 
a single particle might also leave the cluster at some time point. 
In summary, the cluster and particle motion is simulated by Algorithm \ref{alg:mod}. An illustration can be found in figure \ref{sim:movment}

\begin{algorithm}
\begin{algorithmic}
	\State	$\mathbf{Input}$ Initial particle position $p_0=0\in\R^d$ and initial cluster positions $c_{0,i}$, $i \in [N]$ uniformly sampled in $[-b,b]^d$ for some $b>0$. Furthermore $r$: cluster radius, $m_c,\sigma_c^2$: cluster motion parameters, $\sigma_{pc}^2$: variance of the small particle motion inside a cluster, $m_{out}, \sigma_{out}^2$:  particle motion parameters outside  clusters
	
	\For{$t = 0,\ldots,T$}
	\State
	$
	r_t = \min\limits_{n\in[N]} \|p_t-c_{n,t}\|
	$
	\State
	$n_t = \argmin\limits_{n\in[N]} \|p_t-c_{n,t}\|$
	
	\If{$r_t \leq r$}
	\State $\sigma = \sigma_{pc}$
	\State $m = 0$
	\Else
	\State $\sigma = \sigma_{out}$
	\State $m = m_{out}$
	\State $n_t = \mathsf{None}$
	\EndIf
	\State $p_{t+1} = p_t + z_t$, where $z_t \sim \NO(m, \sigma^2 )$
	\State $c_{n,t+1} = c_{n,t} + z_{n,t}$, where $z_{n,t} \sim \NO(m_c, \sigma^2_c )$
	\If{$\hat{n} \neq \mathsf{None}$}
	\State $p_{t+1} = p_{t+1} + (c_{n_t,t+1}-c_{n_t,t})$
	\EndIf
	\EndFor
\end{algorithmic}
\caption{Simulation of particle and cluster movement}\label{alg:mod}
\end{algorithm}

\begin{algorithm}
\caption{EM algorithms for determining the GMM parameters in \eqref{eq:GMM}.}\label{alg:EM}
\begin{algorithmic}
\State Let $f(x, \mu, \sigma) = \frac{1}{2\pi\sigma}\exp\left(-\frac{1}{2}\frac{(x-\mu)^2}{\sigma}\right)$ be the likelyhood of a multivariate Gaussian distribution and 
\State $l(\alpha^r, \mu^r, \sigma^r) = \sum_{t=1}^{T-1}\log\left(\sum_{j\in{in,out}}\alpha^r_j f(\Delta p(t), \mu^r_j, \sigma^r_j)\right)$ be the log likelihood of the GMM for $\{\Delta p(1),\ldots,\Delta p(T-1)\}$.
\State	$\mathbf{Input}$ $ (\Delta p(1),\ldots,\Delta p(T-1))$, Initial values $\alpha^0=(\alpha_{in}^0,\alpha_{in}^0)$, $\mu^0=(\mu^0_{in},\mu^0_{out})$ and $\sigma^r_=(\sigma^0_{in}, \sigma_{out}^0)$. Stopping tolerance $tol$
	\While{$l(\alpha^{r},\mu^r,\sigma^r)-l(\alpha^{r-1},\mu^{r-1},\sigma^{r-1})> tol$}
	\State $\mathbf{E}$ step: \For{$j\in \{in, \ out\}$}
	 $
		\beta_{i,j}^r = \frac{\alpha_j^rf(x_i, \mu_j^r, \sigma_j^r)}{\underset{j\in\{in, out\}}{\sum}\alpha_j^rf(x_i, \mu_j^r, \sigma_j^r)}
	$
	\EndFor
	\State $\mathbf{M}$ step
	\For{$j\in\{in, \ out\}$}
	\begin{align*}
		\alpha_j^{r+1} &= \frac{1}{n} \sum_{i=1}^{n}\beta_{i,j}^r;
		&\mu_j^{r+1}, \sigma_j^{r+1} = \underset{\mu_j, \sigma_j}{\text{argmax}} \sum_{i=1}^n\beta_{i,j}^r\log\left(f(x_i, \mu_j,\sigma_j)\right)
	\end{align*}
	\EndFor

	\EndWhile
\end{algorithmic}
\end{algorithm}

\subsection{GMM for particle motion}
Using the above simulation, the single particle motion in each time step follows a GMM with two modes, i.e., 
\begin{equation} \label{eq:GMM}
	\Delta p(t) \sim	\alpha_{out} \NO(m_{out},\sigma^2_{out}) + \alpha_{in} \NO(m_{in}, \sigma^2_{in})
\end{equation}
where $\alpha_{in}, \alpha_{out} \ge 0$ with $\alpha_{in}+\alpha_{out}=1$ denote the probabilities that 
a particle is inside or outside a cluster at time $t$.
Given measured samples of $p_t$, $t \in [T]$ from a certain number of experiments, we can estimate the
parameters $\alpha_\iota, \mu_\iota, \sigma_\iota$, $\iota\in\{in, out\}$ by the EM algorithm \cite{dempster1977maximum} 
outlined in Algorithm \ref{alg:EM}.

\section{Cluster membership models}\label{sec:clusters}
Given a series of (simulated) cluster and particle positions
$c_{n,t}$ and $p_t$, $n \in [N]$, $t \in [T]$,
the main task is to determine whether a particle lies within a cluster at a certain time $t$.
The problem is even more difficult, since we usually do not know the clusters size $r$, 
but have usually only an upper bound $R$.

In this section, we intend to estimate a probability vector $e \coloneqq (e_t)_{t\in[T]} \in [0,1]^T$
for a particle to be inside some cluster at time $t$ by solving a minimization problem of the form
\[
	\argmin_{e\in [0,1]^T} \ \beta_r L_{r}(e) + \beta_s L_{s}(e), \qquad \beta_r,\beta_s \ge 0.
\]
Here $L_{r}$ takes care of the distance of the particle to the closest cluster
and
$L_{s}$ depends on the local similarity of the movement of the single particle to those of the closest cluster.
In the next subsections, we will propose suitable choices for $L_{r}$ and $L_{s}$. We construct $L_r$ and $L_s$ such that they are convex in $e$ which enables us to use simple gradient descent optimization methods. More specifically we use the adaptive moment estimation optimizer ADAM \cite{kingma2014adam}. Since in our numerical experiments we obtained convergence there was no need for more sophisticated optimization algorithms.

\subsection{Radius function $L_{r}$}
In  general, we do not know the radius $r$ of a cluster,
but can estimate an upper bound such that $r \le R$.  
Let  $\mathcal C_t\coloneqq \{c_{1,t},\ldots,c_{N,t}\}$ denote the set of position of all clusters at time  $t$. 
We define an rough estimate for $e_t$ by
\[
	\bar{e}_{t} \coloneqq
	\begin{cases}
				1 \text{  if  } \text{dist}(p_t,\mathcal C_t) \leq R,\\
				0 \text{  otherwise}.
			\end{cases}
\]
Then we set
\[
	L_{r}(e) \coloneqq \sum_{t=1}^T (e_t-\bar{e}_{t})^2 F \left(\text{dist}(p_t,\mathcal C_t) \right),
\]
where $F$ is the following function which encodes the confidence in our estimation $\bar{e}$:
\begin{align}\label{func:F}
	F(s) \coloneqq \frac{f(s-R)}{f(s-R) + f\left(1-(s-R) \right)}
\end{align}
and
\[
	f(s) \coloneqq 
	\begin{cases}
		e^{-\frac{1}{s}}\text{  for  } s>0,\\
		0 \text{  otherwise}.
	\end{cases}
\]

The function $F$ has the following desirable properties:
\begin{itemize}
\item[i)] $F(s) = 0$ for $s \le R$ where we don't know if the particle is inside a cluster and $F(s) = 1$ for $s \ge R+1$ where we can be certain the particle not inside a cluster,
\item[ii)] $F$ is continuous and monotone increasing.
\end{itemize}
The function is plotted in Figure \ref{fig:plotFG} for $R=0.7$. Since $(e_t-\bar{e}_t)^2$ is convex in $e_t$ and $F \left(\text{dist}(p_t,\mathcal C_t) \right)$ is a constant also $L_r(e)$ is convex in $e$.

\subsection{Similarity functions $L_{s}$}
%
Next we propose two different kinds of local similarity functions.
To this end, we consider trajectories of particles as empirical probability measures. 
More precisely, we define for fixed $k \in \mathbb N$, the sliding window measure of the particles at time $t \in [T]$ by
\begin{equation}\label{local_p}
	\mu_{p,t} \coloneqq \frac{1}{k_t + K_t + 1} \sum_{i=k_t}^{K_t}\delta_{p_{t+i+1}-p_{t+i}},
\end{equation}
where $k_t \coloneqq \max\{-k,1-t\}$ and $K_t \coloneqq \min\{k-1,T-t-1\}$.
In our numerical experiments $k$ will be chosen as $k=6$.
Then $\mu_{p,t}$ can be considered as kind of local particle movement.
Further, we suppose that the parameters $m_\iota, \sigma_\iota$, $\iota \in \{in,out\}$ of the GMM \eqref{eq:GMM}
were estimated by the EM algorithm.

Let $\text{d}: \mathcal P(\mathbb R^d) \times \mathcal P(\mathbb R^d) \to \mathbb R_{\ge 0}$ 
denote a distance function between probability measures.
In this paper, 
we will compare settings using the Wasserstein-2 distance 
and the MMD distance which are explained in Appendix \ref{sec:measures}.

We model the similarity of $\mu_{t}$ to the movement within clusters in two different ways.
The first one takes the GMM parameters $m_{in}, \sigma_{in}$   into account,
while the second one uses a sliding window measure of the nearest cluster centers.

\subsubsection{Similarity between $\mu_{p,t}$ and $\NO(\mu_{in},\sigma^2_{in})$}\label{sec:normal}
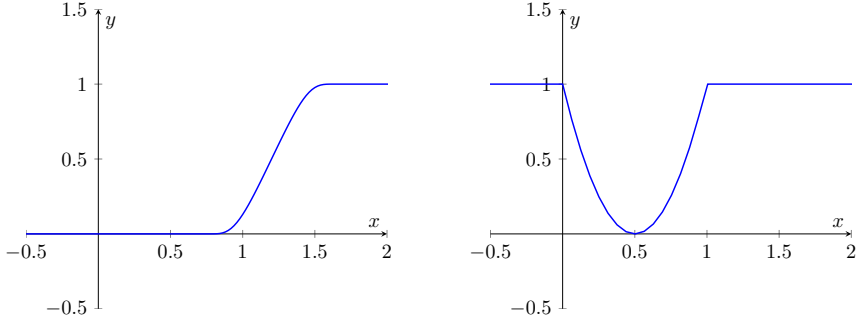
\begin{figure}
\begin{minipage}[t]{0.4\linewidth}
\begin{tikzpicture}[scale=0.7,
  declare function={
        func(\x)= (\x <= 0.0) * (0) + (\x > 0.0) * (exp((-1)/(\x))) ;
        fu(\x) = ((func(\x))/(func(\x)+func(1-(\x))));}]
        \begin{axis}[axis x line=middle, axis y line=middle,
                    ymin=-0.5, ymax=1.5,  ylabel=$y$,
                    xmin=-0.5, xmax=2, xlabel=$x$,samples=501]
                   \addplot [blue,thick] {fu(x-0.7)};
       \end{axis}
\end{tikzpicture}
\end{minipage}
\begin{minipage}[t]{0.3\linewidth}
\begin{tikzpicture}[scale=0.7,
  declare function={
      func(\x)= (x < 0) * (1) + and(\x >= 0, \x < 1) * (4*x^2-4*x+1) + (\x >= 1) * (1); }]
   \begin{axis}[
        axis x line=middle, axis y line=middle,
        ymin=-0.5, ymax=1.5,  ylabel=$y$,
        xmin=-0.5, xmax=2, xlabel=$x$,
        domain=-pi:pi,samples=101]
        \addplot [blue,thick] {func(x)};
   \end{axis}
\end{tikzpicture}
\end{minipage}\caption{Left: Function $F$ from \eqref{func:F} for $R=0.7$. Right: Function $G$ from \eqref{G} for a threshold of $h=0.5$}\label{fig:plotFG}
\end{figure}

For a threshold $h$ which we will determine below, 
we introduce a rough estimate of $e_t$ by
\begin{align}\label{h:gmm}
	\hat{e}_{t} \coloneqq \begin{cases}
		1 \ \text{if} \ \text{d}(\mu_{p,t},\NO(m_{in},\sigma^2_{in} )) \leq h,\\
		0 \ \text{otherwise}.
	\end{cases}
\end{align}
Then we define
\begin{equation} \label{Ls1}
	L_{s}(e) \coloneqq \sum_{t=k+1}^{T-k}(e_t-\hat{e}_{t})^2 G\left(\text{d} \left(\mu_{p,t},\NO(m_{in},\sigma_{in}) \right)\right),
\end{equation}
where 
\begin{equation} \label{G}
	G(s) \coloneqq \min\{1, g(s)\}
\end{equation}
and
$$g(s) \coloneqq \frac 1 h^2 s^2 - \frac 2 h s + 1.$$
Here $g$ is a weight function that encodes the confidence of the estimates $\hat{e}_{t}$. 
The confidence should be small if $\text{d}(\mu_{p,t},\NO(m_{in},\sigma_{in}))$ is close to the threshold. This is true for $g$ since we have, see also figure \ref{fig:plotFG},
\begin{itemize}
\item[i)] $0\leq g(s)\leq 1$ for all $s\in \R$
\item[i)] $g(s) = 0$ for $s = h$ 
\item[ii)] $g(0)=1$ i.e. if $\text{d}(\mu_{p,t},\NO(m_{in},\sigma_{in}))=0$ which means the single particle moves exactly like the clusters
\item[iii)] $g(2h)=1$ i.e. we are confident that the particle is not inside a cluster when its movement differs significantly from that of clusters.
\end{itemize}

Next, we propose a choice for the threshold $h$.
For independent random vectors $X^{\iota}_{-k},\ldots,X^{\iota}_{k-1}$, $\iota \in \{in,out\}$
distributed as $\NO(m_{\iota},\sigma^2_{\iota})$, set 
$\mathbf{X}^\iota \coloneqq (X^{\iota}_{-k},\ldots,X^{\iota}_{k-1})$, and define random measures by 
$$
\boldsymbol{\mu}_{\mathbf{X}^\iota} 
\coloneqq 
\frac{1}{2k}\sum_{i=-k}^{k-1}\delta_{X_i^\iota}, \qquad \iota \in \{in,out\}. 
$$ 
Let the random variables $X^{\iota}_{l}$ be defined on a measure space $\Omega, \Sigma, \alpha$. Then $\Omega:\omega\to\text{d}\left(\boldsymbol{\mu}_{\mathbf{X}^\iota},\NO(m_{in},\sigma^2_{in})\right)$ is a random variable.
With
\begin{align}
	E_{\iota} &\coloneqq 
	\mathbb{E}_{\mathbf{X}^{\iota}}\left[\text{d}\left(\boldsymbol{\mu}_{\mathbf{X}^\iota},\NO(m_{in},\sigma^2_{in})\right) \right]= \int_{\Omega}\text{d}\left(\boldsymbol{\mu}_{\mathbf{X}^\iota},\NO(m_{in},\sigma^2_{in})\right)\dd\alpha,
\end{align}
we suggest as threshold
\[ h \coloneqq \frac{E_{in}+E_{out}}{2}.\]
In the numerical experiments we always use the empirical expected value and variance.

\subsubsection{Similarity between $\mu_{p,t}$ and closest cluster motion}\label{sec:closest}
We can also directly compare the local particle movement to the movement of the closest cluster.
To this end, let 
$$n_t \coloneqq \argmin_{n \in [N]} \ \|c_{n,t}-p_t\|$$ 
be the index of the closest cluster at time point $t$. 
For a window $(t-k,\cdots, t+k-1),$ 
let $\hat n_t$ be the smallest index which occurs most in 
$(n_{t-k},\ldots,n_{t+k-1})$. Then we define the empirical measure
\[
	\mu_{c,t} \coloneqq \frac{1}{k_t + K_t + 1}\sum_{i =k_t}^{K_t}\delta_{c_{\hat n_t, t+i+1}-c_{\hat n_t, t+i}},
\]
where $k_t$ and $K_t$ are given as in \eqref{local_p}.
This describes the movement of the closest cluster at times $(t-k,\ldots, t-k+1)$.
For a threshold $h$ which we will determine below, we introduce a rough estimate of $e_t$ by
\[
	\hat{e}_{t} \coloneqq \begin{cases}
		1 \ \text{if} \ \text{d}(\mu_{p,t},\mu_{c,t})\leq h,\\
		0 \ \text{else}.
	\end{cases}
\]
Then we define
\[
	L_{s}(e) \coloneqq \sum_{t=k+1}^{T-k}(e_t-\hat{e}_{t})^2 G\left( \text{d}(\mu_{p,t}, \mu_{c,t}) \right),
\]
where the confidence function $G$ can be chosen as in \eqref{G}.

Finally, we suggest a threshold $h$.
We assume that we have computed an estimate  $m_c,\sigma_c$ for the cluster movement, 
e.g. by a simple maximum likelihood estimation. 
Let $X_{-k},\ldots,X_{k-1}$ be independent random vectors distributed as $\NO(0,\sigma^2_{in}-\sigma^2_c)$
and  
$Y_{-k},\ldots,Y_{k-1}$ 
be independent random vectors distributed as  $\NO(m_c,\sigma^2_c)$. 
Set $Z_i \coloneqq Y_i + X_i$, $i=-k, \ldots,k-1$, i.e. $Z_i$ models a single particle motion that belongs to this cluster. 
Then we consider 
\begin{align}
	E_{in} &\coloneqq \mathbb{E}_{\mathbf{X},\mathbf{Y}} \left[\text{d}\left(\boldsymbol{\mu}_{\mathbf{Y}}, \boldsymbol{\mu}_{\mathbf{Z}}\right)\right].
\end{align}
Furthermore, let $W_{-k},\ldots,W_{k-1}$ be distributed as $\NO(m_{out},\sigma^2_{out})$ and
\begin{align}
	E_{out} &\coloneqq 
	\mathbb{E}_{\mathbf{X},\mathbf{W}}
	\left[
	\text{d} \left(\boldsymbol{\mu}_{\mathbf{Y}}, \boldsymbol{\mu}_{\mathbf{W}} \right) 
	\right],.
\end{align}
Similarly as before, we define
\[
	h \coloneqq \frac{E_{in}  + E_{out} }{2}.
\]

\begin{remark}
	We used the estimates for the GMM in order to determine the threshold $h$. 
	If we are in a situation, where a GMM is not available, we can also use a clustering method on the set 
	$\{\text{d}(\mu_{p,t},\mu_{c,t}): t\in\{1,\ldots,T\}\}$ to determine $h$.
\end{remark}

\section{Numerical results} \label{sec:num}
All simulations where done in the Julia programming language. For the Wasserstein metric we used the PythonOT for Julia package.
\subsection{EM algorithm}
Table \ref{tab:rmse_results} gives the result of the EM algorithm on simulated data for variuos configurations. It shows that bigger differences in the simulation parameters $\mu_{in}, \mu_{out}$ and $\sigma_{in}, \sigma_{out}$ improves the estimation. If we increase the time step the estimation also improves. The key question is if the estimations suffice in order to get good estimates whether a molecule is inside a cluster in the next section.

\begin{table*}[t]
\centering
	\setlength{\tabcolsep}{3pt}
\begin{tabular}{l| c c| c c || c c | } 
\hline
	Simulation Parameters  &$\sigma_{in}$ & $\sigma_{out}$ & $\mu_{in}$& $\mu_{out}$&$\alpha_{in}$ & $\alpha_{out}$ \\ 
\hline
	\multicolumn{7}{c|}{T: $1000$, r: $0.7$, runs: $10$}\\
 \hline
 
	Simulation parameters & 0.5 & 0.7 & (0.3,0.3)  & 0 & 0.54 & 0.46 \\
	EM estimated 	      & 0.48&0.71&(0.32, 0.29)&(-0.01,-0.01)&0.54&0.46\\
 \hline

	Simulation parameters & 0.5 & 1 & (0.3,0.3)  & 0 & 0.58 & 0.42 \\
	EM estimated 	      &0.50 &1.00&(0.31, 0.30)&(-0.02,-0.04)&0.58&0.42\\
\hline
	Simulation parameters & 0.5 & 0.7 & 0  & 0 & 0.56 & 0.44 \\
	EM estimated 	      &0.45 &0.69&(-0.02, -0.04)&(0.24,-0.09)&0.55&0.45\\
\hline

	\multicolumn{7}{c|}{T: $5000$, r: $0.7$, runs: $10$}\\
 \hline
 
	Simulation parameters & 0.5 & 0.7 & (0.3,0.3)  & 0 & 0.55 & 0.45 \\
	EM estimated 	      & 0.49&0.69&(0.31, 0.30)&(0.02,0.01)&0.52&0.48\\
 \hline

	Simulation parameters & 0.5 & 0.7 & 0 & 0 & 0.57 & 0.34\\
	EM estimated 	      &0.51 &0.73&(0.00,0.00)&(-0.01,0.01)&0.61&0.39\\
\hline
 \end{tabular}
\caption{EM estimations of $\sigma$ and $\mu$}
\label{tab:rmse_results}
 \end{table*}

\subsection{Cluster membership}
For every run of our experiments we simulated data according to algorithm \ref{alg:mod} and estimated the parameters of the corresponding GMM with the EM algorithm. Using the cluster membership models from above we obtain estimates $e=(e_{k+1},\ldots,e_{T-k})$. Note that for $t\leq k, t>T-k$ the measure $\mu_{p,t}$ has fewer than $2k$ points why we omitted them in our analysis. Then let $\hat{e}=(\hat{e}_{k+1},\ldots,\hat{e}_{T-k})$ for
\[
\hat{e}_t= \begin{cases}
0 \text{ if } \ e_t<\frac 1 2\\
1 \text{ else }
\end{cases}
\]
We view $\hat{e}$ as the solution of a classification problem and report the accuracy i.e. the ratio of correctly classified time points. We show the mean accuracy over the runs in table \ref{tab:accuracy}. Detail about the distance functions on $\mathcal{P}_1(\R^d)$ that we used can be found in appendix \ref{sec:measures}.
\subsubsection*{Using the similarity between $\mu_{p,t}$ and $\NO(\mu_{in},\sigma^2_{in})$} 
The first column of table \ref{tab:accuracy} displays the accuracy of the estimation whether the particle is inside a cluster at a given time point. If we only use the radius term $L_r$, i.e., $\beta_r=0$ the peformance is much worse than using only the term $L_s$. Combining $L_r$ and $Ls$ is beneficial for all choices of $L_s$. Generally all choices of $L_s$ perform similarly with the Wasserstein distance based beeing sligthly better than the others. All method perform better when the movement pattern of the single molecule differ stronger from the ones of clusters.
\begin{table*}[t]
\centering
\setlength{\tabcolsep}{3pt}
\begin{tabular}{c c| c c c|| c c c | } 
\hline
	\multicolumn{8}{c|}{$T=1000$, $r=0.7$, $R=1.2$, runs: $20$, $k=6$}\\
 \multicolumn{8}{c|}{$\sigma_{in}=0.5$, $\sigma_{out}=0.7$, $\mu_{in}=(0.3,0.3)$, $\mu_{out}=0$}\\
 \hline
 &  & \multicolumn{3}{|c||}{$\mu_{p,t}$ and $\mathcal{N}(\mu_{in},\sigma^2_{in})$} &\multicolumn{3}{c|}{$\mu_{p,t}$ and $\mu_{c,t}$}\\
&  & \multicolumn{3}{|c||}{Mean Accuracy} &\multicolumn{3}{c|}{Mean Accuracy}\\ 
$\beta_r$ & $\beta_s$&MV  & WS& MMD&MV  & WS& MMD \\

\hline
 $0$&$1$     & 0.77 &0.77&0.74& 0.73 &0.90 &0.90\\
 \hline

 $1$&$1$   &  0.84 & 0.89&0.87  & 0.81 &0.90 &0.94\\
 \hline
 $1$&$0$   &  0.70 & 0.70&0.68  & 0.70 &0.70 &0.69\\
 \hline
 \multicolumn{8}{c|}{$\sigma_{in}=0.5$, $\sigma_{out}=1.0$, $\mu_{in}=(0.3,0.3)$, $\mu_{out}=0$}\\
 \hline
 $0$&$1$   &  0.83 & 0.86&0.79  & 0.72 &0.90 &0.92\\
 \hline
  $1$&$1$   &  0.87 & 0.92&0.90  & 0.82 &0.91 &0.93\\
 \hline
  $1$&$0$   &  0.68 & 0.67&0.68  & 0.68 &0.69 &0.70\\
 \hline
\end{tabular}
\caption{Accuracy of the cluster membership models. Second column: membership model using the similarity between $\mu_{p,t}$ and $\NO(\mu_{in},\sigma^2_{in})$, see Section \ref{sec:normal}. Third column: Accuracy of the cluster membership model using the similarity between $\mu_{p,t}$ and the closest cluster motion, see Section \ref{sec:closest}. MV: Variance and Expected Value distance, WS: Wasserstein-2 Distance, see appendix \ref{sec:measures}}\label{tab:accuracy}
\end{table*}

\subsubsection*{Using the similarity between $\mu_{p,t}$ and closest cluster motion.} 
As demonstrated in the third row of table \ref{tab:accuracy} using the movement 
of the closest cluster has a better accuracy than comparing with $\NO(\mu_{in},\sigma^2_{in})$ for all methods except for the Mean and Variance method. This is especially evindent when we compare the performance without the radius function $L_r$. For the Wasserstein and MMD based method adding $L_r$ does not improve the accuracy by much possibly because the room for improvement is low. Note that there is an inherent error that we have to expect a time point $t$ when a particle enters a cluster but is outside for all time points $t-l$ for $l\leq k$. This could also explain why there is no improvement in accuracy when changing the movement of the single particle from $\NO(0,0.7)$ to $\NO(0,1.0)$.

\section{Conclusions}\label{sec:conc}
We studied how simulated data from a two mode GMM can be used to
to determine if a particle belongs to a certain cluster which itself moves according to a Brownian motion.
Numerical examples give evidence that the approach is suited to solve the task.
Our model can be generalized at several parts. 
The cluster radius may vary between different clusters and maybe also differ over the time.
It may be interesting to consider other distributions instead of the Gaussian ones, 
as e.g. heavy tailed distributions. For an EM algorithm for Student-t distributions, see,
e.g. \cite{HHLS20}. Since our cluster membership model that exploits the similarity of $\mu_{p,t}$ and $\mu_{c,t}$ can be also used by determining the threshold via clustering even settings where the distributions are not known can be tackled. 

Our method shares some similarities with changepoint detection models which use e.g. Wasserstein distance or MMD in that it is a sliding window approach that compares local behavior by distances on the space of probability measures, see e.g. \cite{faber2021watch,wei2022online,li2019scan}. However they compare a window before a certain time point with a window after that time point whereas we compare a window around a time point to either a fixed distribution or to the local behavior of a cluster.
A question is, if there is a practical relevance for such examinations.
Clearly, we want to apply our cluster membership models to 
real-word measurements of tracked molecules and clusters as indicated in the motivation for this work.
The setup for such measurements as well as a reliable molecule position tracking algorithm
to test our model is however a hard task on its own.

\subsection*{Acknowledgement}
Funding by the DFG within the SFB “Tomography Across the Scales” (STE 571/19-1, project
number: 495365311) is gratefully acknowledged. 
Many thanks to J. Hertrich for fruitful discussions.

\bibliographystyle{abbrv}
\bibliography{references}

\appendix
\section{Distances between probability measures}\label{sec:measures}

\subsubsection*{Wasserstein Distance.} 
A popular distance that emerged from optimal transport is the Wasserstein-$p$ distance for $p\geq1$ which is defined as follows see e.g. \cite{villani2009optimal} chapter 6. Let $\mu, \nu\in\P_p(\R^d)$ measures of finite $p$ moment i.e. $\int_{\R^d}\|x\|^p\dd\mu(x)<\infty$. Let $\Gamma(\mu, \nu)\subset \mathcal{P}(\R^d\times \R^d)$ be the set of $\alpha\in\mathcal{P}(\R^d\times \R^d)$ s.t. $\pi^1_{\sharp}\alpha = \mu$ and $\pi^2_{\sharp}\alpha = \nu$ where $\pi^i$ are the projections onto the $i$ component. Then the Wasserstein-$p$ distance is defined as 

\[
	\mathcal{W}_p(\mu,\nu) = \underset{\alpha\in\Gamma(\mu,\nu)}{\mathsf{min}} \ \left(\int_{\R^d\times \R^d}\|x-y\|^p\mathsf{d}\alpha(x,y)\right)^{\frac{1}{p}}.
\]
$\mathcal{W}_p$ is a metric on the space of probability measures $\mathcal{P}_p(X)$ with finite $p$-th moment. Note that finite empirical measures which are the only ones relevant for our application have finite $p$ moment for all $p\geq 1$.

\subsubsection*{Maximum Mean Discrepancy (MMD).} 
Another important distance on $\P(\R^d)$ is the MMD \cite{gretton2006kernel}.
Let $K:\R^d\times \R^d\to \mathbb{R}$ be a kernel. Then for two probability measures $\mu, \nu$ we define
\begin{align*}
	\mathsf{MMD}^2(\mu, \nu) =& \int_{\R^d\times \R^d}K(x,y) \dd\mu(x)\dd\mu(y) - 2\int_{R^d\times \R^d}K(x,y)\dd\mu(x)\dd\nu(y)\\
				  &+ \int_{\R^d\times \R^d}K(x,y)\dd\nu(x)\dd\nu(y).
\end{align*}
We will use as kernel the Riesz kernel $K(x,y) = -\|x-y\|$ for which the MMD defines a metric on $\P_1(\R^d)$.

\subsubsection*{Mean-Variance (MV)} 
Finally, we have also used a simple similarity function based on the variance and expected value
of probability measures.
For $\mu, \nu$, we compare 
\begin{align*}
	\mathbb{E}_{\mu}, \mathbb{E}_{\nu} \\
	\mathsf{cov}(\mu),\mathsf{cov}(\nu)
\end{align*}
where $\mathbb{E}_{\mu}=\int_{\R^d}x\dd \mu\in\R^d$ and $\mathsf{cov}(\mu)=\int_{\R^d}(x-\mathbb{E}_{\mu})(x-\mathbb{E}_{\mu})^T\dd\mu$.
If we assume the covariance matrix to be $\sigma\mathsf{Id}_d$ for $\sigma\in\R$ we also can look at 
\[
	\sqrt{\text{det}(\text{cov}(\mu))}.
\]
We define 
\begin{align}\label{dist:mv}
	d_{m}(\mu,\nu) &= \|\mathbb{E}_{\mu}[x] - \mathbb{E}_{\nu}[x] \|_2  \\
	d_{v}(\mu,\nu) &= |\sqrt{\mathsf{det}(\mathsf{cov}(\mu))} - \sqrt{\mathsf{det}(\mathsf{cov}(\nu))}|.
\end{align}	
	Here we have two functions $L_{s}$ comparing single molecule and cluster movement denoted by $L_{s,mean}$ and $L_{s,var}$ coming from $d_{m}$ and $d_v$, see \eqref{dist:mv}. 
Suppose we estimated $\hat{\mu}_{in}, \hat{\sigma}_{in}, \hat{\mu}_{out}, \hat{\sigma}_{out}$ via EM. 
We then define the threshold $h$ in \eqref{h:gmm} for $L_{s,mean}$ by $h =\frac{\|\hat{\mu}_{in}-\hat{\mu}_{out}\|_2}{2}$. For $L_{s,var}$ we define $h$ by 
$h = \frac{\hat{\sigma}^2_{in} +\hat{\sigma}^2_{out}}{2}$. We then define the $L_s$ in the mean-variance case by 
\[L_s = L_{s,mean} + L_{s,var}.\]

\end{document}